\documentclass[12pt]{article}
\usepackage{amsmath, amssymb, enumerate}

\newtheorem{theorem}{Theorem}[section]
\newtheorem{lemma}[theorem]{Lemma}
\newtheorem{corollary}[theorem]{Corollary} 
\newtheorem{remark}[theorem]{Remark}
\newtheorem{proposition}[theorem]{Proposition}
\newtheorem{assumption}[theorem]{Assumption}

\numberwithin{equation}{section}

\def\grad{{\nabla}}
\def\proof{{\medskip\noindent {\bf Proof. }}}
\def\qed{{\hfill $\square$ \bigskip}}
\def\longproof#1{{\medskip\noindent {\bf Proof #1.}}}

  \def\sC {{\cal C}}
  \def\sF {{\cal F}}
  
  \def\sL {{\cal L}}
 \def\sN {{\cal N}}

 \def\bE {{\Bbb E}}

 \def\bN {{\Bbb N}} 
\def\bP {{\Bbb P}} \def\bQ {{\Bbb Q}} \def\bR {{\Bbb R}}

\begin{document}

\title{Harmonic functions for a class of integro-differential operators.}

\author{Mohammud Foondun}
\date{}
\maketitle

\abstract{
We consider the operator $\sL$ defined on $C^2(\bR^d)$ functions by
\begin{eqnarray*}
\sL f(x)&=&\frac{1}{2}\sum_{i,j=1}^d a_{ij}(x)\frac{\partial^2f(x)}{\partial x_i\partial x_j}+\sum_{i=1}^d b_i(x)\frac{\partial f(x)}{\partial x_i}\\
&+&\int_{\bR^d\backslash\{0\}}[f(x+h)-f(x)-1_{(|h|\leq1)}h\cdot \grad f(x)]n(x,h)dh.
\end{eqnarray*}
\\
Under the assumption that the local part of the operator is uniformly elliptic and with suitable conditions on $n(x,h)$, we establish a Harnack inequality for functions that are nonnegative in $\bR^d$ and harmonic in a domain.  We also show that the Harnack inequality can fail without suitable conditions on $n(x,h)$. A regularity theorem for those nonnegative harmonic functions is also proved.}

\vglue1.0truein
{\small

\noindent{\it Author's address:} Department of Mathematics, The University of Utah,
		155 S. 1400 E. Salt Lake City, UT 84112--0090, USA. {\it Email:}  mohammud@math.utah.edu}

\newpage
%INTRODUCTION
\section{Introduction}
Researchers are increasingly using integro-differential operators (or equivalently, processes with jumps) to model problems from economics and the natural sciences.  For instance, geometric Brownian motion is a standard model for a stock price.  But this model is sometimes not satisfactory because it does not take into account sudden shifts of the stock price.  To model this, one would like to use a process with some jumps, to represent the stock price.  So understanding the properties of those operators is very important.

The purpose of this paper is to consider functions that are harmonic with respect to  the operator $\sL$, where
\begin{eqnarray}\label{def:op}
\sL f(x)&=&\frac{1}{2}\sum_{i,j=1}^d a_{ij}(x)\frac{\partial^2f(x)}{\partial x_i\partial x_j}+\sum_{i=1}^d b_i(x)\frac{\partial f(x)}{\partial x_i} \nonumber\\
&+&\int_{\bR^d}[f(x+h)-f(x)-1_{(|h|\leq1)}h\cdot \grad f(x)]n(x,h)dh
\end{eqnarray}
is defined on $C^2(\bR^d)$ functions.  This is a typical example of a non-local operator, in the sense that the behavior of the harmonic function at a point depends on values of the function at points some distance away rather than just at nearby points.  In probabilistic terms, the local part of $\sL$ corresponds to the continuous part of the process while the non-local part controls the jumps of the process. The jump kernel $n(x,h)$ represents the intensity of jumps from a point $x$ to the point $x+h$ and will be assumed to be nonnegative.

We prove a Harnack inequality as well as a regularity theorem for harmonic functions with respect to the operator $\sL$ without assuming any continuity of the coefficients $a_{ij}$, $b_i$ and of the kernel $n(x,h)$.  We say that a function $u$ is harmonic with respect to $\sL$ in a domain $D$ if $\sL u=0$ in $D$; we give a precise definition in Section 2. Roughly speaking, the Harnack inequality states that the values of a non-negative harmonic function are comparable in a region. In other words, for all $x$ and $y$ lying away from the boundary of $D$, there exists a constant $C$ not depending on $u$ such that 
$$u(x)\leq Cu(y).$$
We show, with the aid of an example, that if $n(x,h)$ does not satisfy some suitable conditions, then a Harnack inequality fails, while under mild conditions on $n(x,h)$, a Harnack inequality holds.

Since the fundamental work of Moser on Harnack inequalities for second order elliptic\cite{M1} and parabolic\cite{M2} partial differential equations with bounded and measurable coefficients, these inequalities have become  increasingly important.  Major contributions to this area have also been made by  Krylov-Safonov\cite{KS1} and Fabes-Stroock\cite{FS}.  While there has been a lot of research on Harnack inequalities for functions that are harmonic with respect to differential operators, not much have been done for non-local operators.  It is only recently that these results have been obtained for harmonic functions associated with purely non-local operators; see \cite{BL1}, \cite{BK1} and \cite{CK}.  The techniques we use to prove the Harnack inequality in this paper are similar to those in \cite{BL1} and in \cite{BK1}, but have their roots in \cite{KS1} where a non-divergence form elliptic operator was considered.  

As for the regularity theorem, we show that there exist $\alpha \in (0,1)$ and a positive constant $C$ not depending on $u$ such that for all $x$ and $y$ lying away from the boundary of $D$, the following holds
\[ |u(x)-u(y)|\leq C\|u\|_\infty |x-y|^\alpha.\]

Continuity estimates of the above type have a long history.  Morrey\cite{Mo} proved such an estimate for second order elliptic partial differential operators in divergence form with bounded coefficients. His result, which was proved in two dimensions only, was independently extended to higher dimensions by DeGiorgi\cite{De} and Nash\cite{N}. Another proof was later given by Moser\cite{M1}.  The corresponding result for operators in non-divergence form was established  by Krylov-Safonov\cite{KS1}.  In \cite{BL2} and \cite{BK2}, the authors considered purely non-local operators and proved a regularity theorem using probabilistic methods. 
It is also interesting to compare our result with the one obtained by ~Mikulevicius-Pragarauskas \cite{MP}.  They considered a parabolic integro-differential operators and obtained a continuity estimate.  However, their result, when specialized to the elliptic case, is a bit weaker than our regularity theorem. In that paper, the jump kernel $n(x,h)$ satisfies a stronger condition than in our paper.  Moreover, our techniques are different.

Another paper which is related to our work here is that of Song-Vondracek $\cite{SV05}$.  Their result is a Harnack inequality for some discontinuous process.  However, the jump kernel considered there is that of a $\alpha$-stable process. Our result thus holds for a much wider class of processes. Related work also include a Harnack inequality for subordinate Brownian motion which has been obtained in \cite{RRV}. %more details and precise reference needed here

The local part of our operator $\sL$ is of non-divergence form. In a forthcoming paper \cite{Fo2}, we consider an operator whose local part is of divergence form and whose jump kernel is symmetric.  In that paper, the problem will be framed in terms of Dirichlet forms and a Harnack inequality together with a regularity theorem will be given.

After stating the results in Section 2, we prove some preliminary estimates in Section 3. In Section 4, we prove a support theorem which is essential to our method. The proof of the Harnack inequality and regularity theorem are given in Section 5 and 6 respectively.  In Section 7, we show that if the jump kernel $n(x,h)$ does not satisfy some suitable conditions, then the Harnack inequality fails.

\section{Statement of results}
We begin this section with some notations and preliminaries.  We use $B(x,r)$ for the open ball of radius $r$ with center $x$.  We also use $|\cdot|$ for the Euclidean norm of points in $\bR^d$, for the norm of vectors and for the norm of matrices.  The letter $c$ with subscripts will denote positive finite constants whose exact values are unimportant.  The Lebesgue measure of a Borel set $A$ will be denoted by $|A|$.

We consider the operator $\sL$ defined by (1.1) and make the following assumptions:

\begin{assumption}\label{assump1}
We assume that the diffusion part of the operator is symmetric and uniformly elliptic and that the ${b_i}s$ are uniformly bounded.  In other words, there exist positive constants $\Lambda_1$ and $\Lambda_2 $ such that
\begin{enumerate}[(a)]
\item the diffusion coefficients $a_{ij}$ satisfy the following
\[\Lambda_1|y|^2 \leq \sum_{i,j=1}^dy_ia_{ij}(x)y_j\leq \Lambda_1^{-1}|y|^2, \hskip 15mm y\in \bR^d, x\in \bR^d,\]
\item \[\sup_i\|b\|_\infty \leq \Lambda_2.\]
\end{enumerate}
\end{assumption}
We let $\sN(\Lambda_1, \Lambda_2)$ denote the set of operators of the form (1.1) satisfying Assumption 2.1.

Besides  nonnegativity, the following assumptions will also be imposed on $n(x,h)$.

\begin{assumption} \label{assump2}
\
\begin{enumerate}[(a)]
\item  There exists a positive constant $K$ such that 
\[ \int_{\bR^d}(|h|^2\wedge 1)n(x,h)dh\leq K, \hskip 10mm \forall x\in\bR^d.\]
\item For any $r\in(0,1]$, any $x_0\in\bR^d$, any $x,\,y\in B(x_0,r/2)$ and $z\in B(x_0,r)^c$, we have $n(x,z-x)\leq k_rn(y,z-y)$, where $k_r$ satisfies $1<k_r\leq kr^{-\beta}$ with $k$ and $\beta$ being positive constants.\\
\end{enumerate}
\end{assumption}

$n(x,h)$ can be thought of as the intensity of the number of jumps from $x$ to $x+h$.  $n(x,z-x)$ thus represents the intensity of the number of jumps from $x$ to $z$.  So Assumption 2.2(b) says that the probability of jumping to a point $z$ is comparable if $x,y$ are relatively far from $z$ but relatively close to  each other.  In Section 7, we show that such an assumption is needed for the Harnack inequality to hold.  

Since our method is probabilistic, we need to work with the Markov process associated with $\sL$. Let $\Omega=D([0,\infty))$ denote the set of paths that are right continuous with left limits, endowed with the Skorokhod topology.  Let $X_t(\omega)=\omega(t)$ for $\omega \in \Omega$ and $\sF_t$ be the right continuous filtration generated by the process $X$. We say a strong Markov process $(\bP^x,X_t)$ is associated with $\sL$ if for each $x$, we have $\bP^x(X_0=x)=1$ and for each $x$ and for each $u\in C^2$ that is bounded and with bounded first and second partial derivatives, $u(X_t)-u(X_0)-\int_0^t\sL u(X_s)ds$ is a local martingale under $\bP^x$.  This is equivalent to saying that $\bP^x$ solves the martingale problem for $\sL$ started at $x$.

We assume that the martingale problem is well posed.  In other words, we assume that the $a_is$ and $b_is$ are continuous so that there exists a unique solution to the martingale problem. We make sure that none of our estimates are dependent on the modulus of continuity of the $a_is$ and $b_is$ so that one can then use an approximation procedure to remove the continuity assumptions.

For any Borel set $A$, let
\[ T_A=\inf \{t: X_t \in A\}, \hskip 20mm \tau_A=\inf \{t: X_t\notin A\},\]
be the first hitting time and first exit time, respectively, of $A$.
We say that the function $u$ is harmonic in a domain $D$ if $u(X_{t\wedge \tau_D})$ is a $\bP^x$-martingale for each $x \in D$.  If $u$ satisfies some regularity conditions and $\sL u=0$ in $D$, it is easy to see that $u$ is harmonic in $D$.
Since our operator contains a non-local part, our process will be have discontinuities.  We write
\[ X_{t-}=\lim_{s\uparrow t} X_s, \hskip 20mm \Delta X_t= X_t-X_{t-}.\]

Our first result  concerns the continuity of harmonic functions. Note that our hypotheses do not require Assumption 2.2(b) to hold.\\
\begin{theorem}\label{theo1}
Suppose	 Assumptions 2.1 and 2.2 (a) hold.  Let $z_0\in \bR^d$ and $R\in(0,1]$. Suppose $u$ is a function which is  bounded in $\bR^d$ and harmonic in $B(z_0,R)$ with respect to $\sL$.  Then there exist $\alpha\in(0,1),\,\,C>0$ depending only on the ${\Lambda_i}'s$ and $K$ such that
\[ |u(x)-u(y)|\leq C\|u\|_\infty \left(\frac{|x-y|}{R}\right)^\alpha, \hskip 10mm x,y\in B(z_0,R/2).\]
\end{theorem}

Our main result is the following Harnack inequality. 
\begin{theorem}\label{theo2}
Suppose Assumptions 2.1 and 2.2 hold.  Let $z_0\in \bR^d$ and $R\in(0,1]$. Suppose $u$ is nonnegative and bounded on $\bR^d$ and harmonic in $B(z_0,R)$ with respect to $\sL$. Then there exists a positive constant $C$ depending on the ${\Lambda_i}'s$, $k$, $\beta$, $R$ and K but not on $z_0$, $u$, or $\|u\|_\infty$ such that
\[ u(x)\leq Cu(y),\hskip 20mm x,y\in B(z_0,R/2).\]
\end{theorem}

\begin{remark}
A chaining argument shows that both results above hold if $R>1$ with $C=C(R)$ depending on R. Theorem 2.3 does not hold for $R>1$ with a constant which is independent of $R$.
\end{remark}

\begin{remark}
For the Harnack inequality, it is essential that $u$ be nonnegative everywhere.  Kassmann~\cite{kas1} has shown that a Harnack inequality can fail for functions $u$ that are harmonic with respect to symmetric stable processes of index $\alpha$ and where $u$ fails to be nonnegative everywhere.
\end{remark}
%Some estimates
\section{Some Estimates}
We start off this section with a proposition which allows us to assume $\Lambda_2=0$ when necessary.                                             The proof is very similar to that of Theorem VI 1.2 in \cite{Ba3}.  See also \cite{St}.  Define
\begin{eqnarray*}
\tilde \sL f(x)&=&\frac{1}{2}\sum_{i,j=1}^d a_{ij}(x)\frac{\partial^2f(x)}{\partial x_i\partial x_j}\\
&+&\int_{\bR^d}[f(x+h)-f(x)-1_{(|h|\leq1)}h\cdot \grad f(x)]n(x,h)dh.
\end{eqnarray*}
\begin{proposition}\label{equivalent}
Suppose $\sL\in\sN(\Lambda_1,\Lambda_2)$.  If there exists a solution, say $\tilde \bP$, to the martingale problem for $\tilde \sL$  started at $x$ where  $\tilde \sL$ is defined as above, then there exists a solution $\bP$ to the martingale problem  for $\sL$ started at $x$.\\
\end{proposition}
%We note that the proof uses Girsanov's theorem and the measure $\bP$ is defined by setting  the restriction of $d \bP/ d\tilde \bP$ to $\sF_t$ to be equal to:
%\[M_t = \exp\left(\int_0^t(b\sigma^{-1}(X_s)dW_s+\frac{1}{2}\int_0^t|b\sigma^{-1}(X_s)|^2ds\right).\]
%where $W_t$ is brownian motion under the measure $\tilde \bP$.  

\begin{proposition}\label{upperbound}
There exist constants $c_1$ and $c_2$ not depending on $x_0$ such that if $r\leq1$, then $\bP^{x}(\tau_{B(x_0,r)}\leq c_1t)\leq tr^{-2}$ for $x\in B(x_0,r)$ and hence
\[\bP^{x}(\tau_{B(x_0,r)}\leq c_2r^2)\leq \frac{1}{2}.\]
\end{proposition}

\proof{Let $u$ be a nonnegative $C^2$ function that is equal to $|x-x_0|^2$ for $|x-x_0|\leq\frac{r}{2}$, which is equal to $r^2$ for $|x-x_0|\geq r$ and such that its first and second derivatives are bounded by $cr$ and $c$ respectively.
Since $\bP^{x}$ solves the martingale problem, we have 
\begin{equation}\label{u:eq1}
\bE^{x}u(X_{t\wedge \tau_{B(x_0,r)}})-u(x_0)=\bE^{x}\int_0^{t\wedge\tau_{B(x_0,r)}}\sL u(X_s)ds.  
\end{equation}
\\
Let us write the operator $\sL$ as $\sL=\sL_c+\sL_d$ where
\[\sL_cu(x)=\frac{1}{2}\sum_{i,j=1}^d a_{ij}(x)\frac{\partial^2u(x)}{\partial x_i\partial x_j}+\sum_{i=1}^d b_i(x)\frac{\partial u(x)}{\partial x_i},\] 
and 
\[\sL_du(x)=\int[u(x+h)-u(x)-1_{(|h|\leq1)}h\cdot \grad u(x)]n(x,h)dh.\]
Since the first and second derivatives of $u(x)$ are bounded, we have
$\sL_cu(x)\leq c_3$ for $x\in B(x_0,r)$ and hence
\begin{equation}\label{u:eq2}
\Big|\int_0^{t\wedge\tau_{B(x_0,r)}}\sL_cu(X_s)ds\Big|\leq c_3t.\\
\end{equation}
Now let us look at $\sL_du(x)$ for $x\in B(x_0,r)$
\begin{eqnarray*}
|\sL_du(x)|&=&\Big|\int[u(x+h)-u(x)-1_{|h|\leq1}h\cdot\grad u(x)]n(x,h)dh\Big|\\
 &\leq&\Big|\int_{|h|\leq 1}[u(x+h)-u(x)-1_{(|h|\leq1)}h\cdot\grad u(x)]n(x,h)dh\Big|\\
 %&+&\Big|\int_{r\leq |h| \leq 1}[u(x+h)-u(x)-1_{(|h|\leq1)}h\cdot\grad u(x)]n(x,h)dh\Big|\\
 &+&\Big|\int_{|h|\geq1}[u(x+h)-u(x)]n(x,h)dh\Big|\\
 &=& I_1+I_2.
 \end{eqnarray*}
\\
By our assumptions and the fact that the second derivatives of $u(x)$ are bounded, we get
\begin{eqnarray*}
I_1&=&\Big|\int_{|h|\leq 1}[u(x+h)-u(x)-1_{(|h|\leq1)}h\cdot\grad u(x)]n(x,h)dh\Big|\\
&\leq&c_3\int_{|h|\leq 1}|h|^2  \| D^2u\|_\infty n(x,h)dh\leq c_4,\\
%I_2&=&\Big|\int_{r\leq |h| \leq 1}[u(x+h)-u(x)-1_{(|h|\leq1)}h\cdot\grad u(x)]n(x,h)dh\Big|\\
%& \leq &\int_{r\leq |h| \leq 1}2\|u\|_\infty n(x,h)dh+\int_{r\leq |h| \leq 1}|h\cdot \grad u(x)|n(x,h)dh\\
%&\leq& \int_{r\leq |h| \leq 1}2r^2n(x,h)dh+c\int_{r\leq |h| \leq 1}r|h|n(x,h)dh\\
%&\leq &c_5+c_6\int_{r\leq |h| \leq 1}|h|^2n(x,h)dh\\
%& \leq &c_7,\\
I_2&=&\Big|\int_{|h|\geq1}[u(x+h)-u(x)]n(x,h)dh\Big|\\
&\leq&2\|u\|_\infty \int_{|h|\geq 1}n(x,h)dh \leq c_5.\\
\end{eqnarray*}
Hence we have $\left|\int_0^{t\wedge \tau_{B(x_0,r)}}\sL_du(X_s)ds\right|\leq c_6t$.
This together with \eqref{u:eq1} and \eqref{u:eq2} yield
\[\bE^{x}u(X_{t\wedge \tau_{B(x_0,r)}})\leq c_7t,\]
and so from $r^2\bP^{x}(\tau_{B(x_0,r)}\leq t)\leq \bE^xu(X_{t\wedge \tau_{B(x_0,r)}})$, we get the first part of the proposition.  The second part is obtained by choosing $t=\frac{1}{2}r^2$.  
\qed

%The following proposition concerns scaling.  The proof is straightforward and is omitted here.
%\\
%\begin{proposition}
%Let $(\bP^x,X_t)$ be a solution of the martingale problem with operator $\sL$, let $Y_t=aX_{t/{a^2}}$.  Define $\bQ^x=\bP^{x/a}$. Then $(\bQ^x,Y_t)$ is a solution to the martingale problem with operator $\tilde{\sL}$ where
%\begin{eqnarray*}
%\tilde{\sL}f(x)&=&\frac{1}{2}\sum_{i,j=1}^d a_{ij}(\frac{x}{a})\frac{\partial^2f(x)}{\partial x_i\partial x_j}+\sum_{i=1}^d \frac{1}{a}b_i(\frac{x}{a})\frac{\partial f(x)}{\partial x_i}\\
%&+&\int_{\bR^d}[f(x+h)-f(x)-1_{(|h|\leq a)}h\cdot \grad f(x)]\tilde{n}(x,h)dh,\\
%\end{eqnarray*}
%and $\tilde{n}(x,h)=\frac{1}{a^{2+d}}n(\frac{x}{a},\frac{h}{a})$.

%\end{proposition}

We have the following L\'evy system formula:
\begin{proposition}\label{L-S}
If $A$ and $B$ are disjoint Borel sets, then for each $x$,
\begin{equation}\label{L-S}
\sum_{s\leq t}1_{(X_{s-}\in A,X_s\in B)}-\int_0^t\int_B1_A(X_s)n(X_s,u-X_s)duds
\end{equation}
is a $\bP^x$-martingale.
\end{proposition}
The proof is identical to that of the purely non-local operator and can be found in \cite{BL1}.
\\
\begin{lemma}\label{Lem: exp}
There exist $c_1$ and $c_2$ such that if $r\leq\frac{1}{2}$,
\begin{enumerate}[(a)]
\item $\bE^{x}\tau_{B({x_0},r)}\geq c_1r^2$ for $x\in B(x_0,r/2)$ and 
 \item $\bE^{x}\tau_{B({x_0},r)}\leq c_2r^2$ for $x\in B(x_0,r)$.
\end{enumerate}
\end{lemma}

\proof{ By Proposition \ref{upperbound}, there exists $c_3$ such that $\bP^{x}(\tau_{B({x_0},r)}\leq c_3r^2)\leq\bP^x(\tau_{B(x,r/2)}\leq c_3r^2)\leq \frac{1}{2}$.  The first inequality follows by writing
\begin{equation*}
\bE^x\tau_{B({x_0},r)}\geq c_3r^2\bP^{x}(\tau_{B({x},r)}\geq c_3r^2),
\end{equation*}
and using the above. Now let us look at the proof of the second inequality. For simplicity we assume that $\sL \in \sN(\Lambda_1,0)$.  The general case follows by using Proposition \ref{equivalent} (and a change of measure argument). Since $\bP^{x}$ solves the martingale problem, we have 
\begin{equation}\label{e:eq 1}
\bE^{x}u(X_{t\wedge \tau_{B(x_0,r)}})-u(x_0)=\bE^{x}\int_0^{t\wedge\tau_{B(x_0,r)}}\sL u(X_s)ds.
\end{equation}\\
As before, let us write $\sL=\sL_c+\sL_d$.  Let us choose a bounded smooth function $u(x)$ so that  $\begin{displaystyle}u(x)=|x-x_0|^2\end{displaystyle}$ for $x\in B(x_0,2)$ and $u(x)$ equals some constant greater than $4$ outside the ball $B(x_0,4)$. Some calculus shows that $\sum_{i,j=1}^d\partial_{ij}u(x)=\sum_{i=1}^d\partial_{ii}u(x)$ and is a constant for $x\in B(x_0,2)$.  This and the uniform ellipticity of the local part of $\sL$ implies that there exists a positive constant $c_4$ such that $\sL_cu(X_s)\geq c_4$ whenever $X_s\in B(x_0,r)$.

To deal with the non-local part, we write
\begin{eqnarray*}
 \sL_du(x)&=& \int_{|h|\leq1}[u(x+h)-u(x)-h\cdot \grad u(x)]n(x,h)dh\\
 &+&\int_{|h|>1}[u(x+h)-u(x)]n(x,h)dh\\
&=& I_1+I_2.
\end{eqnarray*}
Note that for $|h|\leq1$, we have $x+h \in B(x_0,3/2)$ for $x\in B(x_0,r)$, so by convexity and the fact that $n(x,h)\geq 0$, we obtain $I_1\geq 0$.  As for the second term, we have
\begin{eqnarray*}
I_2&\geq& \int_{|h|>1}[|x+h-x_0|^2-|x-x_0|^2]n(x,h)dh\\
&\geq&0.
\end{eqnarray*} 
The facts that $x\in B(x_0,r)$ and $|h|>1$ imply that $x+h\notin B(x_0,r)$ which means that the integrand is always non-negative.
Combining the above, we have $\sL_du(X_s)\geq 0$ whenever $X_s\in B(x_0,r)$ and hence
\[\bE^{x}\int_0^{t\wedge\tau_{B(x_0,r)}}\sL u(X_s)ds \geq c_4\bE^{x}(t\wedge \tau_{B(x_0,r)}).\]
\\
Now the left hand side of (\ref{e:eq 1}) satisfies:
\[\bE^{x}u(X_{t\wedge \tau_{B(x_0,r)}})-u(x_0) \leq r^2.\]
\\
Combining the above and letting $t\rightarrow \infty$, we get the second inequality.
\qed
\begin{corollary}\label{Cor:exp}
For each $p\geq 1$, there exists $c_1$ depending on $p$ such that for $r\leq\frac{1}{2}$ and $x \in B(x_0,r)$,
\[\bE^{x}(\tau_{B(x_0,r)}^p)\leq c_1r^{2p}.\]
\end{corollary}

\proof{Note that
\begin{equation*}
\bE^{x}(\tau_{B(x_0,r)})\geq t\bP^{x}(\tau_{B(x_0,r)}\geq t).
\end{equation*}
Letting $t=2r^2$ and using Lemma \ref{Lem: exp}, we obtain $\bP^{x}(\tau_{B(x_0,r)}\geq 2r^2)\leq \frac{1}{2}$}. If $\theta_t$ is the shift operator from the theory of Markov processes, then by the Markov property
\begin{eqnarray*}
\bP^{x}(\tau_{B(x_0,r)}\geq (m+1)r^2)&\leq& \bP^{x}(\tau_{B(x_0,r)} \geq mr^2,\tau_{B(x_0,r)}\circ \theta_{mr^2}\geq r^2)\\
&=&\bE^{x}[\bP^{X_{mr^2}}(\tau_{B(x_0,r)}\geq r^2); \tau_{B(x_0,r)}\geq mr^2]\\
&\leq&\frac{1}{2}\bP^{x}(\tau_{B(x_0,r)}\geq mr^2).
\end{eqnarray*}
By induction $\bP^{x}(\tau_{B(x_0,r)}\geq mr^2)\leq (\frac{1}{2})^m$. The required result then follows easily from this.}
\qed

%\begin{theorem}\label{theo: lepelMar}
%Consider  $\sL \in \sN(\Lambda_1, \Lambda_2)$ and suppose that the jump kernel $n(x,h)$ satisfies the following: $\begin{displaystyle} \int_{\bR^d}[|h|^21_{(|h|\leq1)}+|h|1_{(|h|>1)}]n(x,h)dh < M \end{displaystyle}$, for some constant $M$.  Let $x\in \bR^d$.  If there exists a solution $\bP^x$ to the martingale problem associated with $\sL \in \sN(\Lambda_1, \Lambda_2)$, then for any bounded measurable function $f$, the following holds:
%\begin{equation}\label{kry}
%\bE^x \int_0^{t\wedge \tau_{B(x_0,r)}} f(X_s) ds\leq N\|f\|_{L^d(B(x_0,r))}
%\end{equation}
%where $N$ depends on the $\Lambda_is$, $t$, $r$ and $M$.\\
%\end{theorem}

The following result was first proved in the continuous case by Krylov in \cite{kry1}.  Since then, this inequality has been extended for diffusions with jumps. See for instance Theorem III$_{15}$ of \cite{LP}. What follows is a consequence of Corollary 2 of \cite{MP}.

\begin{proposition}\label{Cor: lepelMar}
Consider  $\sL \in \sN(\Lambda_1, \Lambda_2)$ and suppose that Assumption 2.2(a) holds.  Let $x\in \bR^d$ and $R\in (0,1]$.  If  $(\bP^x,X_t)$ is a solution to the martingale problem associated with $\sL \in \sN(\Lambda_1, \Lambda_2)$, then for any bounded measurable function $f$, the following holds:
\begin{equation}\label{kry1}
\bE^x \int_0^{\tau_{B(x_0,R)}} f(X_s) ds\leq NR \|f\|_{L^d(B(x_0,R)),}
\end{equation}
where $N$ depends on the $\Lambda_i's$, and $K$.
\end{proposition}

%\proof{
%Let $n_0(x,h)=n(x,h)1_{(|h|\leq 2r)}$.  Also let $\sL_0f(x)$ be defined similarly to (1.1), but with $n_0(x,h)$ (instead of $n(x,h)$) as jump kernel. If $(\overline{\bP}^x,\overline{X}_t)$ is the solution of the martingale problem associated with $\sL_0$, then by Theorem 3.6, we have
%\begin{equation}\label{kry2}
 %\overline{\bE}^x \int_0^{t\wedge \tau_{B(x_0,r)}} f(\overline{X}_s) ds\leq N\|f\|_{L^d(B(x_0,r))}.\\
%\end{equation}
%Since $n_0(x,h)\leq n(x,h)$, and Assumption 2.2(a) gives
%\[ \int_{\bR^d}(n(x,h)-n_0(x,h))dh \leq c, \hskip 10mm {\rm for\,\,\,all}\,\,\,x\in \bR^d,\]
%we can now use Meyer's construction to recover the process $X_t$ such that $(\bP^x,X_t)$ will be a solution to the martingale problem associated with the operator $\sL$.
%Notice that the integral on the left hand side of (3.6) is obtained by summing $f$ along the trajectories of $X_t$ inside the ball $B(x_0,r)$ only.  Since $\overline{X}_t$ differs from $X_t$ in that $X_t$ has a jump of size larger than $2r$ (in which case, $X_t$ will be outside $B(x_0,r)$), inequality (3.6) holds for $(\bP^x, X_t)$ also. 
%Finally, a closer look at the proof of Theorem \ref{theo: lepelMar} reveals the dependence of $N$ (from (\ref{kry})) on the radius of the ball $B(x_0,r)$.
%}
%qed

\begin{remark}
In the above, the constant $N$ depends on the non-local part only through the constant $K$.  It also does not depend on the radius $R$ but the upper bound does depend on the radius as shown in \eqref{kry1}.
\end{remark}

The following gives a lower bound on the probability that our process hits a set $A$ before leaving a ball.  However, the set $A$ should fill most of the ball.  We will extend this result in the next section.
\begin{proposition}\label{prop:hit}
Let $\epsilon>0$ and $r\in(0,1/2]$. If $x\in B(x_0,\frac{r}{2})$, $A\subset B(x_0,r)$ and $|B(x_0,r)-A|\leq \epsilon$, then $\bP^{x}(T_A\leq \tau_{B(x_0,r)})\geq \rho(r,\epsilon)$, where there exists some $\epsilon_0$ such that $\rho(r,\epsilon)>0$ for $0<\epsilon\leq \epsilon_0$.
\end{proposition}
\proof{ From inequality (\ref{kry1}), we have 
\begin{eqnarray*}
\left|\bE^{x}\int_0^{\tau_{B(x_0,r)}}1_{A^c}(X_s)ds\right|&=&\left|\bE^{x}\int_0^{\tau_{B(x_0,r)}}1_{(B(x_0,r)-A)}(X_s)ds\right|\\
&\leq&Nr|B(x_0,r)-A|^{1/d}\\
&\leq&Nr\epsilon^{1/d}.
\end{eqnarray*}
So we can write
\begin{eqnarray*}
\bE^{x}\tau_{B(x_0,r )}&\leq&\bE^{x}(\tau_{B(x_0,r )};T_A\leq \tau_{B(x_0,r)})+\bE^{x}\int_0^{\tau_{B(x_0,r)}}1_{A^c}(X_s)ds \nonumber\\
&\leq&\left(\bE^{x}\tau_{B(x_0,r )}^2\right)^{1/2}\left(\bP^{x}(T_A\leq \tau_{B(x_0,r)})\right)^{1/2}+Nr\epsilon^{1/d} .
\end{eqnarray*}
From Lemma \ref{Lem: exp} and Corollary \ref{Cor:exp} we have $\bE^{x}\tau_{B(x_0,r )}\geq c_1r^2$ and $\bE^{x}\tau_{B(x_0,r )}^2\leq c_2r^4$. So the above yields
\[\bP^{x}(T_A\le \tau_{B(x_0,r)})\geq \left(\frac{c_1r^2-Nr\epsilon^{1/d}}{c_3r^2}\right)^2.\]
The proposition  is then proved with $\rho(r,\epsilon)=(\frac{c_1r^2-Nr\epsilon^{1/d}}{c_3r^2})^2.$
\qed

The following will be used only in the proof of the Harnack inequality. So far, this is the only place where we use Assumption 2.2(b).
\begin{proposition}\label{prop: assump2}
Under Assumption 2.2, there exists a constant $c_1$ which depends on $K$, such that if $r\leq 1/2$, $z\in B(x_0,\frac{r}{4})$ and $H$ is a bounded non-negative function supported in $B(x_0,r)^c$, then
\begin{equation}\label{a0}
\bE^{x_0}H(X_{\tau_{B(x_0,\frac{r}{2})}})\leq c_1k_r\bE^zH(X_{\tau_{B(x_0,\frac{r}{2})}}).
\end{equation}
\end{proposition}

\proof{ By linearity and a limit argument, it suffices to consider only $H(x)=1_C(x)$ for a set $C$ contained in $B(x_0,r)^c$.
From Assumption 2.2(b), we have $n(w,v-w)\leq k_rn(y,v-y)$ for all $w,y\,\,\in B(x_0,\frac{r}{2})$ and $v\in B(x_0,r)^c$. Hence we have,
\begin{equation}\label{a1}
 \sup_{y\in B(x_0,\frac {r}{2})}n(y,v-y)\leq k_r \inf_{y\in B(x_0,\frac {r}{2})}n(y,v-y).
\end{equation}
By optional stopping and the L\'evy system formula, we have 
\begin{eqnarray*}
\bE^z1_{(X_{t\wedge \tau_{B(x_0,\frac{r}{2})}}\in C)}&=&\bE^z\sum_{s\leq{t\wedge \tau_{B(x_0,\frac{r}{2})}}}1_{(|X_s-X_{s-}|\geq \frac{r}{2}, X_s\in C)}.\\
&=&\bE^z\int_0^{t\wedge \tau_{B(x_0,\frac{r}{2})}}\int_C n(X_s,v-X_s)dvds.\\
&\geq&\bE^z(t\wedge \tau_{B(x_0,\frac{r}{2})})\int_C \inf_{y\in B(x_0,\frac{r}{2})}n(y,v-y)dv.
\end{eqnarray*}
Letting $t\to\infty$ and using the dominated convergence theorem on the left and monotone convergence on the right, we obtain
\[\bP^z(X_{\tau_{B(x_0,\frac{r}{2})}}\in C)\geq \bE^z\tau_{B(x_0,\frac{r}{2})}\int_C\inf_{y\in{B(x_0,\frac{r}{2})}}n(y,v-y)dv.\]
Since $ \bE^z\tau_{B(x_0,\frac{r}{2})}\geq \bE^z\tau_{B(z,\frac{r}{4})}$, we have 
\begin{equation}\label{a2}
\bP^z(X_{\tau_{B(x_0,\frac{r}{2})}}\in C)\geq \bE^z\tau_{B(z,\frac{r}{4})}\int_C\inf_{y\in{B(x_0,\frac{r}{2})}}n(y,v-y)dv. 
\end{equation}
Similarly we have 
\begin{equation}\label{a3}
\bP^{x_0}(X_{\tau_{B(x_0,\frac{r}{2})}}\in C)\leq \bE^{x_0}\tau_{B(x_0,\frac{r}{2})}\int_C\sup_{y\in{B(x_0,\frac{r}{2})}}n(y,v-y)dv.
\end{equation}
Combining inequalities (\ref{a1}), (\ref{a2}) and (\ref{a3}) and using Lemma \ref{Lem: exp}, we get our result.
\qed

Our process is a discontinuous one consisting of small jumps as well as big jumps.  In many cases it is more convenient to discard the big jumps and add them later.  This can be done by using a construction which is due to Meyer \cite{Me}.  We will use this in the next section for the proof of the support theorem.
\subsection*{Meyer's construction:}
Suppose that we have two jump kernels $n_0(x,h)$ and $n(x,h)$ with $n_0(x,h)\leq n(x,h)$ and such that for all $x\in \bR^d$,
\[N(x)=\int_{\bR^d}(n(x,h)-n_0(x,h))dh \leq c.\]
Let $\sL$ and $\sL_0$ be the operators corresponding to the kernels $n(x,h)$ and $n_0(x,h)$ respectively. If $\overline X^0_t$ is the process corresponding to the operator $\sL_0$, then we can construct a process $X_t$ corresponding to the operator $\sL$ as follows.  Let $S_1$ be an exponential random variable of parameter 1 independent of $X_t$, let $C_t=\int_0^tN(X_s)ds$, and let $U_1$ be the first time that $C_t$ exceeds $S_1$.
At the time $U_1$, we introduce a jump from $X_{U_1-}$ to $y$, where $y$ is chosen at random according to the following distribution:
\[ \frac{n(X_{U_1-},h)-n_0(X_{U_1-},h)}{N(X_{U_1-})}dh.\]
This procedure is repeated using an independent exponential variable $S_2$.  Since $N(x)$ is finite, this procedure adds only a finite number of big jumps on each finite time intervals.  In \cite{Me}, it is proved that the new process corresponds to the operator $\sL$.

%Support theorem
\section{Support theorem} 

The main result of this section is the support theorem.  Before stating and proving this result,   we present some ideas which will be crucial for its proof. More precisely, we will represent the solution of the martingale problem as a solution to a stochastic differential equation. We begin by representing the jumps of our discontinuous process as a function of a Poisson point process.

Suppose that $\bP^x$ is a solution to the martingale problem associated with $\sL \in \sN(\Lambda_1,\Lambda_2)$ started at $x$. Let $Y_s$ be the point process associated with $X_s$, that is, $Y_s=\Delta X_s$ if $\Delta X_s \neq 0$ and $0$ otherwise.  Then there exists a measurable function $F(x,z)$ such that $Y_s=F(X_{s-},\hat{Y}_s)$ where $\hat{Y}_s$ is a Poisson point process with intensity measure $\tilde{\lambda}$.  To make this statement more precise, let $F(x,A)=\{ F(x,z): z\in A \}$ and define
\begin{equation*}
N_t(A)=\sum_{s\leq t}1_{(\Delta X_s\in F(X_{s-}, A))}.
\end{equation*}
Then under $\bP^x$, $N_t(\cdot)$ is a Poisson point process with intensity measure $\tilde{\lambda}$.  Moreover the measure $\tilde{\lambda}$ satisfies the following
\begin{enumerate}[(a)]
\item 
\begin{equation}\label{prop1} 
\int (|z|^2\wedge 1)\,\tilde{\lambda} (dz)<\infty,
 \end{equation}
\item
 \begin{equation}\label{prop2}        
\int 1_A(h)n(x,h)dh=\int 1_A(F(x,z))\tilde{\lambda}(dz).
\end{equation}
\end{enumerate}
\begin{remark}
The second condition above gives the relationship between the jump kernel $n(x,h)$ and the Poisson process $N_t(\cdot)$.  Moreover, the indicator function in (\ref{prop2}) can be replaced by a larger class of functions. For a more precise statement and proof of the above, see Theorem 12 in \cite{EL}. See also Chapter XIV of \cite{Ja}.
\end{remark}
We now relate $(\bP^x, X)$, the solution of the martingale problem to that of a stochastic differential equation. Set $\mu([0,t]\times A)=N_t(A)$ and $\nu([0,t]\times A)=t\tilde{\lambda}(A)$.  Let $W_t$ be a  Brownian motion with respect to the filtration $\sF_t$.  Then $X_t$ solves the following stochastic differential equation
\begin{eqnarray}\label{rep}
dX_t&=&\sigma(X_t)dW_t+b(X_t)dt+\int _{|F(X_{t-},z)|\leq1}F(X_{t-},z)(\mu-\nu)(dz,dt) \nonumber \\
& &+\int_{|F(X_{t-},z)|>1}F(X_{t-},z)\mu(dz,dt), \hskip 20mm X_0=x,
\end{eqnarray}
where $\sigma\sigma^T$ has $a_{ij}$ as entries and $\sigma^T$ denotes the transpose of $\sigma$.  The above has been taken from \cite{LP}. Chapter XIV of \cite{Ja} contains more information about this relation.  In fact, according to Theorem II$_{10}$ of \cite{LP}, this representation holds under a more stringent condition on the big jumps of the process $\hat{Y}_t$ (see Property M in \cite{LP}).  Since the proof of the theorem below involves dealing with small jumps only and then adding the big jumps later, this does not  affect our result (see the proof below).  Here is our support theorem:
\begin{theorem}\label{theo:support} 
Suppose $\sL \in\sN(\Lambda_1,\Lambda_2)$ and $\bP^{x_0}$  is a solution to the martingale problem for $\sL$ started at $x_0$. Let  $\epsilon>0$ and suppose that $\phi:[0,t_0]\to \bR^d$ is differentiable with $\phi(0)=x_0$.  There exist constants $c_1$, $c_2$ and $c_3$ depending on $\Lambda_1$, $\Lambda_2$, $t_0$, $K$ and $\sup_{t\leq t_0}|\phi'(t)|$ but not on $\epsilon$ such that for all $\lambda>0$,
\begin{equation}\label{S0}
 \bP^{x_0}(\sup_{t\leq t_0}|X_t-\phi(t)|< \epsilon)\ge c_1\Big[1-\exp\big[-\lambda \epsilon+\frac{\lambda^2t_0}{2}(c_2+c_3e^{|\lambda|})\big]\Big]^2.
\end{equation}
\end{theorem}

The support theorem says that the graph of $X_s$ stays inside an $\epsilon$-tube  about $\phi$.  In other words, if $G_\phi^\epsilon=\{(s,y): |y-\phi(s)|< \epsilon,\,s\leq t \}$, then $\{(s, X_s): s\leq t \}$ is contained in $G_\phi^\epsilon$ with a positive probability.

\longproof{of theorem \ref{theo:support}} We will assume $\sL\in \sN(\Lambda_1,0)$.  The general case follows from Proposition \ref{equivalent}.  We first consider the special case when the jump kernel of $\sL$ is defined by $n_0(x,h)=n(x,h)1_{(|h|<1)}$ and denote the corresponding process by $\overline{X}_t$.  We will later use Meyer's construction to remove this restriction. We now use the stochastic differential equation representation of the solution to the martingale problem. In other words, we use the fact that $\overline{X}_t$ satisfies (\ref{rep}).  Since our process do not have jumps of size greater than $1$, the last term of ($\ref{rep}$) can be taken to be identically zero.  More precisely, replacing $1_A(h)$ by $|h|1_{(|h|>1)}$ and $n(x,h)$ by $n_0(x,h)$ in (\ref{prop2}), we obtain $\int_{|F(x,z)|>1}|F(x,z)|\tilde{\lambda}(dz)=0$.
Define a new measure $\bQ$ by 
\begin{align}\label{S1}
\frac{d\bQ}{d\bP^{x_0}}=\exp\Big[-\int_0^{t_0}\phi'(s)\sigma^{-1}(\overline{X}_{s-})dW_s-\frac{1}{2}\int_0^{t_0}|\phi'(s)\sigma^{-1}(\overline{X}_{s-})|^2ds\Big].
\end{align}
Let
\begin{equation*}
Z_t=\overline{X}_t-\int_0^t \int_{|F(\overline{X}_{s-},z)|\le1}F(\overline{X}_{s-},z)(\mu-\nu)(dz,ds).
\end{equation*}
We see that
\begin{eqnarray*}
{\left\langle-\int_0^t\phi'(s)\sigma^{-1}(\overline{X}_{s-})dW_s,Z_t \right\rangle} 
& = &\left \langle-\int_0^t\phi'(s)\sigma^{-1}(\overline{X}_{s-})dW_s, \int_0^t \sigma(\overline{X}_{s-})dW_s\right\rangle\\
& = & -\int_0^t\phi'(s)ds = -\phi(t)+\phi(0).
\end{eqnarray*}
\\
So by Girsanov's theorem, under $\bQ$, each component of $Z_t$ is a semi-martingale. If
\begin{eqnarray*}
\widehat{W_t}&=&\int_0^t\sigma^{-1}(\overline{X}_{s-})d\overline{X}_s-\int_0^t\int_{|F(\overline{X}_{s-},z)|\le1}\sigma^{-1}(\overline{X}_{s-})F(\overline{X}_{s-},z)(\mu-\nu)(dz,ds) \nonumber \\
&-&\int_0^t\sigma^{-1}(\overline{X}_{s-})\phi'(s)ds,
\end{eqnarray*}
then $\widehat{W_t}$ is a continuous martingale and $d\langle\widehat{W}_t^i,\widehat{W}_t^j\rangle$=$\delta_{ij}dt$ under $\bQ$. Hence $\widehat{W_t}$ is  a $d$-dimensional Brownian motion under $\bQ$.  Note
\begin{equation}\label{S2}
d(\overline{X}_t-\phi(t))=\sigma(\overline{X}_{t-})d\widehat{W_t}+\int_{|F(\overline{X}_{t-},z)|\le1}F(\overline{X}_{t-},z)(\mu-\nu)(dt,dz).\\
\end{equation}
\\
If we prove the following:
\begin{equation}\label{S3}
\bQ(\sup_{t\leq t_0}|\overline{X}_t-\phi(t)|<\epsilon)\geq c_4,
\end{equation}\\
then the theorem will be proved, for if $A$ is the event $\{ \sup_{s\leq t_0}|\overline{X}_s-\phi(s)|< \epsilon\}$, then
\begin{equation}\label{eq:depend}
c_5\leq \bQ(A)=\int_A(d\bQ/d\bP^{x_0})d\bP^{x_0} \leq(\bE^{x_0}(d\bQ/d\bP^{x_0})^2)^{\frac{1}{2}}( \bP^{x_0}(A))^{\frac{1}{2}}.
\end{equation}
The theorem then follows easily by noting the $d\bQ/d\bP^{x_0}$ has a finite second moment which is bounded by a constant depending on $t_0$, $\Lambda_1$ and $\sup_{t\leq t_0}|\phi'(t)|$; see page 188 of \cite{Ba3}. Now let us look at the proof of (\ref{S3}).  Let us write the left hand side of (\ref{S2}) as $dD_t$ i.e, $D_t:=\overline{X}_t-\phi(t)$.

Let $\lambda$ be a constant to be chosen later. Define
\begin{equation}\label{S4'}
N_t=\lambda D_t-\frac{\lambda^2}{2}\int_0^t|\sigma(\overline{X}_{s-})|^2ds-\int_0^t\int_{|z|\leq1}(e^{\lambda z}-1-\lambda z)n(\overline{X}_{s-},z)dzds.
\end{equation}
Set $Y_t^\lambda=e^{N_t}$.  Then, by Ito's formula (for processes with jumps), we obtain
\begin{eqnarray*}
Y_t^\lambda&=&1+\lambda\int_0^te^{N_{s-}}dD_t-\int_0^t\int_{|z|\leq1}e^{N_{s-}}(e^{\lambda z}-1-\lambda z)n(\overline{X}_{s-},z)dzds\\
& &+\sum_{s\leq t}[e^{N_s}-e^{N_{s-}}-e^{N_s}\Delta N_s].\\
\end{eqnarray*}
 From  Assumptions 2.1(a) and 2.2(a), there exist constants $c_5$ and $c_6$ such that
\begin{equation*}
\int_0^t|\sigma(\overline{X}_s)|^2ds\leq c_6t, 
\end{equation*}
and 
\begin{eqnarray*}
\int_{|z|\leq1}(e^{\lambda z}-1-\lambda z)n(\overline{X}_{s-},z)dz&\leq&\frac{\lambda^2}{2}e^{|\lambda|}\int_{|z|\leq 1}|z|^2n(\overline{X}_{s-},z)dz\\
&\leq&c_7\frac{\lambda^2}{2}e^{|\lambda|}.
\end{eqnarray*}
By noting that $\Delta N_s=\lambda \Delta D_s$, and using Theorem 10 of \cite{LP} together with the above, we see that $Y_t^\lambda$ is a martingale. The above bounds, together with ($\ref{S4'}$) also yield
\begin{equation*}
 \bQ(\sup_{t\leq {t_0}}|D_t|\ge \epsilon)\leq \bQ(\sup_{t\leq {t_0}}e^{N_t}\ge \exp[\lambda \epsilon-\frac{\lambda^2}{2}c_6t_0-\frac{\lambda^2}{2}e^{|\lambda|}c_7t_0]).
\end{equation*}
Since $Y_t^\lambda=e^{N_t}$, we can apply Doob's inequality as follows:
\begin{eqnarray}\label{S4}
 \bQ(\sup_{t\leq {t_0}}|D_t|\ge \epsilon)&\leq&\bQ(\sup_{t\leq t_0}Y^\lambda_t\geq \exp[\lambda \epsilon-\frac{\lambda^2}{2}c_6t_0-\frac{\lambda^2}{2}e^{|\lambda|}c_7t_0])\nonumber \\
&\leq&\bE_{\bQ}Y^\lambda_{t_0}\exp[-\lambda \epsilon+\frac{\lambda^2}{2}c_6t_0+\frac{\lambda^2}{2}e^{|\lambda|}c_7t_0]\nonumber \\
&=&\bE_{\bQ}Y^\lambda_0\exp[-\lambda \epsilon+\frac{\lambda^2}{2}c_6t_0+\frac{\lambda^2}{2}e^{|\lambda|}c_7t_0]\nonumber \\
&\leq&\exp[-\lambda \epsilon+\frac{\lambda^2}{2}c_6t_0+\frac{\lambda^2}{2}e^{|\lambda|}c_7t_0].
\end{eqnarray}
From the above we conclude that
\begin{equation}\label{S5}
 \bQ(\sup_{t\leq {t_0}}|D_t|<\epsilon)\geq 1-\exp\big[-\lambda \epsilon+\frac{\lambda^2t_0}{2}(c_6+c_7e^{|\lambda|})\big].\\
\end{equation}
We now use Meyer's construction to recover the process $X_t$ so that $(\bP^x,X_t)$ is a solution to the martingale problem associated with the operator $\sL$ whose jump kernel satisfies the weaker Assumption 2.2(a).  The trajectories of $X_t$ now have jumps greater than 1. Recall that $U_1$ is the first time that $C_t$ exceeds $S_1$ where $S_1$ is an exponential random variable with parameter 1.  More precisely, we have
\begin{eqnarray}\label{S6}
\bQ(\sup_{t\leq {t_0}}|X_t-\phi(t)|< \epsilon)&=&\bQ(\sup_{t\leq {t_0}}|X_t-\phi(t)|< \epsilon; U_1\leq t_0)+ \bQ(\sup_{t\leq {t_0}}|X_t-\phi(t)|< \epsilon; U_1>t_0)\nonumber \\
 &\geq&\bQ(\sup_{t\leq {t_0}}|\overline{X}_t-\phi(t)|<\epsilon)\bQ(U_1>t_0).
 %&\leq&\bQ(U_1\leq t_0)+c_7\bQ(U_1>t_0).
\end{eqnarray}
%&=& 1-[1-\bQ(U_1\leq t_0)][1-\bQ(\sup_{t\leq t_0}|\overline{X}_t-\phi(t)|\geq \epsilon)].
Using the fact that
\[\bQ(U_1\leq t_0)\leq \bQ(S_1\leq(\sup N)t_0)=1-e^{-(\sup N)t_0},\]
inequality (\ref{S6}) reduces to $\bQ(\sup_{t\leq {t_0}}|X_t-\phi(t)|<\epsilon)\geq c_8\bQ(\sup_{t\leq {t_0}}|\overline{X}_t-\phi(t)|<\epsilon)$ for some positive constant $c_8$.  This inequality, together with \eqref{S5} and \eqref{eq:depend} complete the proof.
\qed}
\begin{remark}\label{rem:support}
By taking $\lambda=\epsilon^2$ in inequality \eqref{S0}, we obtain upon choosing $\epsilon$ small enough,
\begin{equation*}
\bP^{x_0}(\sup_{t\leq t_0}|X_t-\phi(t)|< \epsilon)\ge c_1[1-e^{-\epsilon^3/2}]^2.
\end{equation*}
We now use the fact that $1-e^{-x}\geq (1-e^{-1})x$ whenever $0\leq x \leq1$ to obtain
\begin{equation}\label{eq:rem:support}
\bP^{x_0}(\sup_{t\leq {t_0}}|X_t-\phi(t)|<\epsilon) \geq c_2\epsilon^6,
\end{equation}
for some positive constant $c_2$ not depending on $\epsilon$.  However, $c_2$ does depend on $\phi$  via $\sup_{t\leq t_0}|\phi'(t)|$; see \eqref{eq:depend} and the discussion following it.
\end{remark}
We now present a corollary of the above support theorem. 

\begin{corollary}\label{corol:support} 
Suppose $\sL \in\sN(\Lambda_1,\Lambda_2)$ and $\bP^{x_0}$  is a solution to the martingale problem for $\sL$ started at $x_0$. Let  $\epsilon>0$ and suppose that $\phi:[0,t_0]\to \bR^d$ is  continuous with $\phi(0)=x_0$.  There exist constants $c_1$, $c_2$ and $c_3$ depending on $\epsilon$, $\Lambda_1$, $\Lambda_2$, $t_0$ and the modulus of continuity of $\phi$ such that for all $\lambda>0$,
\begin{equation}\label{cor:S0}
 \bP^{x_0}(\sup_{t\leq t_0}|X_t-\phi(t)|< \epsilon)\ge c_1\Big[1-\exp\big[-\lambda \epsilon+\frac{\lambda^2t_0}{2}(c_2+c_3e^{|\lambda|})\big]\Big]^2.
\end{equation}
\end{corollary}
\proof{ Let us choose a differentiable function $\phi_d$ with derivative bounded by say $c_4$ and such that $\sup_{s\leq t}|\phi(s)-\phi_d(s)|\leq \epsilon$.  Moreover, we can choose $\phi_d(s)$ such that $\|\phi_d'\|_\infty$ depends only on $t$, $\epsilon$ and the modulus of continuity of $\phi$; see Page 60 of \cite{Ba2}.  Hence proving the following
\begin{equation*}
\bP^{x_0}(\sup_{t\leq t_0}|X_t-\phi_d(t)|<\epsilon)\ge c_6\Big[1-\exp\big[-\lambda \epsilon+\frac{\lambda^2t_0}{2}(c_7+c_8e^{|\lambda|})\big]\Big]^2
\end{equation*}
will imply (\ref{cor:S0}) but with $2\epsilon$ instead of $\epsilon$. But the above inequality follows from Theorem \ref{theo:support}. Hence the corollary is proved.
\qed}

\begin{remark}
The above corollary only requires the function $\phi$ to be continuous but the downside of this generalization is that we can longer keep track of the dependence of the constants on $\epsilon$.
\end{remark}

Let $Q(x,r)$ denote the cube of side length $r$ centered at $x$. If $R_i$ denotes a cube with side length $r$, then $\hat{R_i}$ also denotes a cube with the same center but with side length $r/3$. The next result is not a probabilistic result.  It enables us to decompose $Q(0,1)$ into smaller subcubes such that a subset $A$ of $Q(0,1)$ fills a percentage of each of the smaller subcubes. Since this is Proposition V.7.2 of \cite{Ba3}, we do not include a proof here.
\begin{proposition}\label{decomp}
Let $q\in(0,1)$. If $A\subseteq Q(0,1)$ and $|A| \leq q$, then there exists $D$ such that (i) $D$ is the union of cubes $\hat{R}_i$ such that the interiors of the $R_i$ are pairwise disjoint, (ii) $|A|\leq q|D\cap Q(0,1)|$, and (iii) for each i, $|A \cap R_i|\geq q|R_i|$.
\end{proposition}
 A  corollary of the support theorem is the following:
\begin{corollary}\label{Cor:Support}
Let $r\in(0,R)$ and $R\in(0,1]$.  Let $y\in Q(0,R)$ with $dist (y,\partial Q(0,R)) \ge r$, $\sL \in \sN(\Lambda_1,\Lambda_2)$, and $\bP$ be the solution to the martingale problem started at $y$. If $Q(z,r)\subseteq Q(0,R)$, then $\bP(T_{Q(z,r)}\le \tau_{Q(0,R)})\geq \zeta(r)$ where $\zeta(r)>0$ depends only on $r$, $K$ and the $\Lambda_is$.
\end{corollary}
The above two results together with the Proposition \ref{prop:hit} are the main ingredients in obtaining the estimate below. The proof is essentially the same as that of Theorem V7.4 in \cite{Ba3} so we omit it here.
\begin{proposition}\label{Cor:Hit}
There exists a non-decreasing function $\psi:(0,1)\rightarrow(0,1)$ such that if $B\subseteq Q(0,R)$, $|B|> 0$, $R\in(0,1]$ and $x\in Q(0,R/2)$, then
\[ \bP^x(T_B\le \tau_{Q(0,R)})\geq \psi(|B|/R^d).\]
\end{proposition}
We now give a different version of the above proposition.  This will allow us to use balls instead of cubes.
\begin{corollary}\label{cu-ba}
There exists a non-decreasing function $\phi$, such that if $B\subseteq B(0,R)$, $|B|> 0$, $R\in(0,1]$ and $x\in B(0,R/2)$, then
\[ \bP^x(T_B\le \tau_{B(0,R)})\geq \phi(|B|/R^d).\]
\end{corollary}
\proof{
For simplicity, we assume $d=2$.  Higher dimensional cases differ only in notation.  Let $k$ be a large positive integer and let $R_{ij}$ be squares of the form $[(i-1)R/k,iR/k]\times[(j-1)R/k,jR/k]$, where $i,j \in \{-k+1,...,-1,0,1,...,k\}$. 

Take $\epsilon>0$ small enough and $k$ sufficiently large so that 
\[\sC=\{R_{ij}: |R_{ij}\cap B|>0,\,\, {\rm and} \,\,R_{ij}\subset B(0,(1-\epsilon)R)\}\]
is nonempty.  Let $M$ be the number of elements in $\sC$.
Let $R_{ij}^\ast$ be the cube with the same center as $R_{ij}$ but side length half as long. Let $D=\cup_{R_{ij}\in \sC}R_{ij}^\ast$. Pick $z\in R_{ij}^\ast$, where $R_{ij}$ satisfies $|R_{ij} \cap B|\geq \frac{|B|}{M}$. We can choose $k$ larger if necessary so that we can find such a cube. Then using Proposition \ref{Cor:Hit} and the fact that $|R_{ij}|=\frac{R^2}{k^2}$, we have 
\begin{eqnarray*}
\bP^z(T_{R_{ij} \cap B}<\tau_{R_{ij}})&\geq& \psi(|B|k^2/MR^2)\\
&\geq& \psi(|B|/R^2),
\end{eqnarray*}
where the last inequality is obtained by  noting that $M\leq k^2$.  Since dist$(x, \partial B(0,R))\geq R/2$ and $D\subset B(0,R)$, we can use Corollary \ref{Cor:Support} to obtain
\[ \bP^x(T_D<\tau_{B(0,R)})\geq c_1,\]
where $c_1$ is a constant.
Using the Markov property and the above inequalities, we obtain
\begin{eqnarray*}
\bP^x(T_B<\tau_{B(0,R)})&\geq& \bE^x[\bP^{X_{T_{R_{ij}^\ast}}}(T_B<\tau_{R_{ij}});T_{R_{ij}^\ast}<\tau_{B(0,R)}]\\
&\geq&c_1\psi(|B|/R^2).
\end{eqnarray*}
}
\qed

%Proof of the Regularity Theorem
\section{The Regularity Theorem}
Now we are ready to prove the regularity theorem.
\\
\longproof{of Theorem 2.3}
Let us suppose $u$ is bounded by $M$ in $\bR^d$ and $z_1\in B(z_0,R/2)$. Set
\[ r_n= \theta_2 \rho^{n}, \hskip 10mm s_n=\theta_1a^n, \hskip10mm {\rm for\,\,\,\,}n\in \bN,\]
where $a < 1$, $\rho<1/2$, and $\theta_1 \ge 2M$ are constants to be chosen later. We choose $\theta_2$ small enough that $B(z_1, 2r_1)\subset B(z_0,R/2)$. Write $B_n=B(z_1,r_n)$ and $\tau_n=\tau_{B_n}$. Set
\[M_n=\sup_{x\in B_n} u(x), \hskip 10mm m_n=\inf_{x\in B_n} u(x).\]
We will use induction to show that $M_n-m_n \leq s_n$ for all $n$.  The H\"older continuity at $z_1$ follows from this.  Let $n_0$ be a positive number to be chosen later.  Suppose $M_i-m_i\leq s_i$ for all $i=1,2,...,n$, where $n\geq n_0$; we want to show 
\[M_{n+1}-m_{n+1} \leq s_{n+1}.\]
\
Let $\epsilon >0$ and choose $z,\,y\in B_{n+1}$ such that $u(y)\leq m_{n+1}+\epsilon$ and $u(z)\geq M_{n+1}-\epsilon$.  We will show that $u(z)-u(y)\leq s_{n+1}$ and since $\epsilon>0$ is arbritary, this will imply $M_{n+1}-m_{n+1} \leq s_{n+1}$ as desired.

Let
\[ A_n=\{ x \in B_{n}: u(z)\leq(M_n+m_n)/2\}.\]
We may suppose that $|A_n|/|B_{n}|\geq 1/2$, for if not, we can look at $M_n-u$ instead. Let $A$ be  a compact subset of $A_n$ such that $|A|/|B_n|\geq 1/3$. Corollary \ref{cu-ba} gives the following
\begin{equation}\label{R0}
\bP^x(T_A\leq \tau_n)\geq c_1,
\end{equation}
where $c_1$ is a constant and $x\in B_{n+1}$. 
Let $z,\,y\in B_{n+1}$.  By optional stopping,
\begin{eqnarray}\label{R1}
u(z)-u(y)&=&\bE^z[u(X_{T_A})-u(y); T_A \le \tau_n] \nonumber\\
&+& \bE^z[u(X_{\tau_n})-u(y); \tau_n \le T_A, X_{\tau_n}\in B_{n-1}] \nonumber\\
&+& \sum_{i=1}^{n-2} \bE^z[u(X_{\tau_n})-u(y); \tau_n \le T_A, X_{\tau_n}\in B_{n-i-1}-B_{n-i}] \nonumber\\
&+& \bE^z[u(X_{\tau_n})-u(y);  \tau_n \le T_A, X_{\tau_n}\notin B_1]\nonumber\\
&=& I_1+I_2+I_3+I_4.
\end{eqnarray}
By the L\'evy system formula, and Lemma \ref{Lem: exp} (see the proof of Proposition 3.5 of \cite{BL1}) , there exist $c_2$ and $c_3$ such that 
\begin{eqnarray}\label{R2}
\sup_{y\in B_{n+1}}\bP^y(X_{\tau_n}\notin B_{n-i})&\leq&\sup_{y\in B_{n+1}}\bE^y\tau_n\int_{|h|>r_{n-i}-r_{n+1}}n(y,h)\,dh\nonumber\\
&=&\sup_{y\in B_{n+1}}\bE^y\tau_n[\int_{|h|>1}n(y,h)\,dh+\int_{1\geq|h|>r_{n-i}-r_{n+1}}n(y,h)\,dh]\nonumber\\
&\leq& c_2r_n^2+c_3\left( \frac{\rho^i}{1-\rho^i}\right)^2.
\end{eqnarray}
%In particular, we have
%\begin{equation}\label{R3}
%\bP^y(X_{\tau_n}\notin B_1)\leq c_2 r_n^2.\\
%\end{equation}
The first term on the right of (\ref{R1}) is bounded as follows 
\begin{equation}\label{R3}
I_1\leq \left(\frac{M_n+m_n}{2}-m_n\right)\bP^y(T_A \le \tau_n) \leq \frac{1}{2}s_n\bP^y(T_A\le \tau_n).\\
\end{equation}
As for the second term, we have
\begin{equation}\label{R4}
I_2 \leq (M_{n-1}-m_{n-1})\bP^y(\tau_n \le T_A)\leq s_{n-1}(1-\bP^y(T_A \leq \tau_n)).\\
\end{equation}
To bound the third term, we choose $\rho=\frac{\sqrt{a}}{2}\wedge \sqrt{\frac{3c_1a}{128c_3}}$ and note that
\begin{eqnarray*}
\sum_{i=1}^{n-2}s_{n-i-1}&=&s_{n-1}\sum_{i=1}^{n-2}a^{-i}\\
&\leq&s_{n-1}[\frac{a^2}{a^n(1-a)}]
\end{eqnarray*}
and
\begin{equation*}
\sum_{i=1}^{n-2}s_{n-i-1}\rho^{2i}\leq s_{n-1}[\frac{\rho^2/a}{1-\rho^2/a}].
\end{equation*}

Using (\ref{R2}) and the above, the third term is bounded by
\begin{eqnarray*}
 \sum_{i=1}^{n-2}(M_{n-i-1}&-&m_{n-i-1})\bP^y(X_{\tau_n}\notin B_{n-i}) \nonumber\\
 &\leq& c_2r_n^2\sum_{i=1}^{n-2} s_{n-i-1}+c_3\sum_{i=1}^{n-2}s_{n-i-1}\rho^{2i}\nonumber\\
 &\leq&s_{n-1}[\frac{c_2a^2\theta_2^2\rho^{2n}}{a^{n}(1-a)}+\frac{c_3\rho^2/a}{1-\rho^2/a}].\nonumber
 \end{eqnarray*}
By our choice of $\rho$, we obtain $1-\rho^2/a\geq 3/4$,  $\rho^{2n}/a^n\leq 1/2^{2n}$ and $\rho^2/a\leq \frac{3c_1}{128c_3}$ so that the above reduces to
 \begin{eqnarray*}
 \sum_{i=1}^{n-2}(M_{n-i-1}&-&m_{n-i-1})\bP^y(X_{\tau_n}\notin B_{n-i})\nonumber\\
 &\leq&s_{n-1}[\frac{a^2\theta_2^2c_4}{1-a}+\frac{4c_3\rho^2/a}{3}]\nonumber\\
 &\leq&s_{n-1}[\frac{a^2\theta_2^2c_4}{1-a}+\frac{c_1}{32}]\nonumber.
 \end{eqnarray*}
 We also choose $\theta_2$ smaller if necessary so that $\begin{displaystyle}\theta_2\leq \frac{1}{4}\sqrt{\frac{c_1(1-a)}{2a^2c_4}} \end{displaystyle}$ and obtain
\begin{equation}\label{R5}
I_3\leq \frac{s_{n-1}c_1}{16}.
\end{equation}
Using (\ref{R2}) again, we see that the fourth term is bounded by 
\begin{eqnarray*}
2M \bP^y(X_{\tau_n} \notin B_1)&\leq&2M[c_2r_n^2+c_3\rho^{2(n-1)}]\\
 &\leq& \theta_1[c_2a^{4n}\theta_2^2+c_3a^{4n-4}].
 \end{eqnarray*}
 By choosing $n_0$ bigger if necessary and recalling that $a<1$, we obtain for $n\geq n_0$, 
 \begin{equation}\label{R6}
 I_4\leq \frac{s_{n-1}c_1}{8}.
 \end{equation}
Inequalities (\ref{R0})-(\ref{R6}) give the following:
 \begin{eqnarray*}
 u(y)-u(z)&\leq& \frac{1}{2}as_{n-1}\bP^y(T_A\le \tau_n)+ s_{n-1}(1-\bP^y(T_A \leq \tau_n)) + s_{n-1}[\frac{c_1}{16}+\frac{c_1}{8}]. \nonumber\\
 \end{eqnarray*}
Using the fact that $a$ is less than one, we obtain
 \begin{eqnarray*}
 u(z)-u(y)&\leq& \frac{s_n}{a}\Big[1-\frac{\bP^y(T_A<\tau_n)}{2}+\frac{c_1}{16}+\frac{c_1}{8}\Big]\\
 &\leq&\frac{s_n}{a}[1-\frac{5c_1}{16}].\\
 \end{eqnarray*}
We now choose $a$ as follows:
 \[  a=\sqrt{1-\frac{5c_1}{16}}.\]
 This yields
 \begin{equation}
 u(z)-u(y)\leq s_na=s_{n+1}. 
 \end{equation}
The continuity estimate now follows from \cite{M1}.
\qed
 }

\section{ Proof of the Harnack Inequality}

\longproof{ of Theorem 2.4} By looking at $u+\epsilon$ and letting $\epsilon \downarrow 0$. We may suppose that $u$ is bounded below by a positive constant. Also, by looking at $au$, for a suitable $a$, we may suppose that $\inf_{B(z_0,R/2)}u \in [1/4,1]$.
 \
 We want to bound $u$ above in $B(z_0,R/2)$ by a constant not depending on $u$.  Our proof is by contradiction.
 
 Since $u$ is continuous, we can choose $z_1\in B(z_0,R/2)$ such that $u(z_1)=\frac{1}{3}$.  Let $r_i=r_1Ri^{-2}$ where $r_1<\frac{1}{2}$ is a chosen constant so that $\sum_{i=1}r_i< R/8.$ Recall that from  Proposition \ref{prop: assump2}, there exists $c_1$ such that if $r<\frac{1}{2}$, $y\in B(x,r/4)$ and $H$ is a bounded non-negative function supported in $B(x,r)^c$, then
\begin{equation}\label{Har1}
\bE^xH(X_{\tau_{B(x,r/2)}}) \leq c_1k_{r}\bE^yH(X_{\tau_{B(x,r/2)}}).
\end{equation}
For inequality (\ref{Har1}) to hold, we need Assumption 2.2(b).  Let  $\eta$ be a constant to be chosen later. Also, let $\xi$ be a constant defined as follows
\[\xi= \frac{1}{2}\wedge \frac{\eta}{c_1}.\]
Let $c_2$, $c_3$ and $c_4$ be positive constants to be chosen later.  Once these constants have been chosen, we suppose that there exists $x_1\in B(z_0,R/2)$ with $h(x_1)=K_1$ for some $K_1$ large enough so that the following is satisfied: 
\begin{equation}\label{Har2}
\frac{\xi K_1e^{c_2j}r_j^{\beta+6} c_3c_4}{k}\geq 2,
\end{equation}
for all $j$. This is possible because of the fact that $r_j=r_1Rj^{-2}$. The constants $k$  and $\beta$ are taken from Assumption 2.2(b).  

We will show that there exists a sequence $\{(x_j,K_j)\}$ with $x_{j+1}\in \overline{B(x_j,r_j)}\subset B(x_j,2r_j)\subset B(z_0,3R/4)$ with:
\begin{equation}\label{Har3}
K_j=u(x_j)\quad{\rm and}\quad K_j\geq K_1e^{c_2j}.
\end{equation}
This would imply that $K_j\rightarrow \infty$ as $j\rightarrow \infty$ contradicting the fact that $u$ is bounded.  Suppose that we already have $x_1,x_2,...,x_i$ such that (\ref{Har3}) is satisfied.  We will show that there exists $x_{i+1}\in \overline{B(x_i,r_i)}\subset B(x_i,2r_i)$ such that $K_{i+1}=u(x_{i+1})$ and $K_{i+1}\geq K_1e^{c_2(i+1)}$. Then by induction, (\ref{Har3}) will hold for all $j$.  Define
\begin{equation*}
 A=\{ y\in B(x_i,r_i/4);  u(y)\geq \frac{\xi K_ir_i^{\beta}}{k} \}.
\end{equation*}
We are going to show that $\begin{displaystyle} |A|\leq \frac{1}{2}|B(x_i,r_i/4)| \end{displaystyle}$.  To prove this fact, we suppose the contrary.  Choose a compact set $A' \subset A$ with $\begin{displaystyle} |A'| > \frac{1}{2}|B(x_i,r_i/4)| \end{displaystyle}$.  Note that upon choosing $r_1$ smaller if necessary, we can use \eqref{eq:rem:support} to obtain 
\begin{equation*}
\bP^{z_1}(T_{B(x_i,r_i/4)}<\tau_{B(z_0,R)})\geq c_5r_i^6,
\end{equation*}
where $c_5$ is independent of $r_i$.  To see this, consider \eqref{eq:rem:support} with  $\epsilon=r_i/8$ and let $\phi$ be a line segment joining $z_1$ and $x_i$ ($\phi(0)=z_1$ and $\phi(t_0)=x_i$). Since $z_1, x_i\in B(z_0, R)$, $|\phi'(t)|$ is bounded by a constant which is independent of $i$. Hence $c_5$ is also independent of $i$; see \eqref{eq:depend} and the discussion following it.

Hence, using the strong Markov property, we can write
\begin{eqnarray*}
\bP^{z_1}(T_{A'}<\tau_{B(z_0,R)})
&\geq&\bE^{z_1}[\bP^{X_{T_{B(x_i,r_i/4)}}}(T_{A'}<\tau_{B(x_i,r_i)});T_{B(x_i,r_i/4)}<\tau_{B(z_0, R)}]\\
&\geq&\phi\left(\frac{|A'|}{|B(x_i,r_i)|}\right)\bP^{z_1}(T_{B(x_i,r_i/4)}<\tau_{B(z_0,R)})\\
&\geq&\phi((1/2)^{2d+1})c_5r_i^6.
\end{eqnarray*}

 We now take $c_3=\phi\left((1/2)^{2d+1}\right)$ and $c_4=c_5$.
By optional stopping,  the above inequality and the fact that $u(X_{t\wedge T_{A'}})$ is right continuous, we obtain
\begin{eqnarray*}
\frac{1}{3}&=& u(z_1)\geq \bE^{z_1}[u(X_{T_{A'}\wedge \tau_{B(z_0,R)}}); T_{A'} < \tau_{B(z_0,R)}]\\
&\geq&  \frac{\xi K_ir_i^\beta}{k} \bP^{z_1}(T_{A'}<\tau_{B(z_0,R)})\\
&\geq&\frac{\xi K_1e^{c_2i}r_i^{\beta+6} c_3c_4}{k}\\
&\geq& 2.
\end{eqnarray*}
This is a contradiction.  Therefore $\begin{displaystyle}|A|\leq \frac{1}{2} |B(x_i,r_i/4)|\end{displaystyle}$.
So we can find a compact set $E$ such that $E \subset B(x_i,r_i/4)-A$ and $|E|\geq \frac{1}{3}|B(x_i,r_i/4)|$. Let us write $\tau_{r_i}$ for $\tau_{B(x_i,r_i/2)}$.  From Corollary \ref{cu-ba}, we have $\bP^{x_i}(T_E<\tau_{r_i})\geq c_6$ where $c_6$ is some positive constant.

Let $M=\sup_{B(x_i,r_i)}u(x)$.  We then have
\begin{eqnarray}\label{Har4}
K_i=u(x_i)&=&\bE^{x_i}[u(X_{T_E\wedge \tau_{r_i}}); T_E < \tau_{r_i}]\nonumber \\
&+& \bE^{x_i}[u(X_{T_E\wedge \tau_{r_i}}); T_E > \tau_{r_i}, X_{\tau_{r_i}}\in B(x_i,r_i)]\nonumber \\
&+&\bE^{x_i}[u(X_{T_E\wedge \tau_{r_i}}); T_E > \tau_{r_i}, X_{\tau_{r_i}}\notin B(x_i,r_i)]\nonumber \\
&=& I_1+I_2+I_3.
\end{eqnarray}
Writing $p_i=\bP^{x_i}(T_E<\tau_{r_i})$, we see that the first two  terms are easily bounded as follows:
\[ I_1 \leq \frac{\xi K_ip_ir_i^\beta}{k}, \hskip 15mm {\rm and}\hskip 15mmI_2\leq M(1-p_i).\]
To bound the third term, we prove $\bE^{x_i}[u(X_{\tau_{r_i}}); X_{\tau_{r_i}} \notin B(x_i,r_i)] \leq \eta K_i$.  If not, then by using (\ref{Har1}), we will have, for all $y\in B(x_i,r_i/4)$, 
\begin{eqnarray*}
u(y)&\geq& \bE^yu(X_{\tau_{r_i}})\geq \bE^y[u(X_{\tau_{r_i}}); X_{\tau_{r_i}}\notin B(x_i,r_i)]\\
&\geq& \frac{1}{c_1k_{r_i}}\bE^{x_i}[u(X_{\tau_{r_i}});X_{\tau_{r_i}}\notin B(x_i,r_i)]>\frac{\eta K_i}{c_1k_{r_i}}>\frac{\xi K_ir_i^\beta}{k},
\end{eqnarray*}
contradicting the fact that $|A|\leq \frac{1}{2}|B(x_i,r_i/4)|$. Hence
\[ I_3 \leq \eta K_i.\]
So (\ref{Har4}) becomes
\begin{eqnarray*}
K_i\leq \frac{\xi K_ip_ir_i^\beta}{k}+M(1-p_i)+\eta K_i\\
%&\leq& \xi K_ip_i+M(1-p_i)+\eta K_i
\end{eqnarray*}
or 
\begin{eqnarray}\label{Har5}
\frac{M}{K_i}&\geq&\frac{1-\eta-\xi p_ir_i^{\beta}/k}{1-p_i}\nonumber\\
&=&1+\frac{(1-\xi r_i^\beta/k)p_i-\eta}{1-p_i}.
\end{eqnarray}
Choosing $\eta=\frac{c_6}{4}$ and using the definition of $\xi$ together with the fact that $p_i\geq c_6$ and $r_i^\beta/k<1$, we see that  there exists a positive $L$, such that inequality (\ref{Har5}) reduces to $M\geq K_i(1+L)$.  Therefore there exists $x_{i+1}\in \overline{B(x_i,r_i)}$ with $u(x_{i+1}) \geq K_i(1+L).$ Setting $K_{i+1}=u(x_{i+1})$, we see that 
\begin{eqnarray*}
K_{i+1}&\geq&K_i(1+L)\\
&=&K_ie^{\log(1+L)}.
\end{eqnarray*}
The condition \eqref{Har3} is thus satisfied provided we choose $c_2=\log(1+L)$.  Finally, note that the fact that $\sum_{i=1}r_i< \frac{R}{8} $ implies that $B(x_i,2r_i)\subset B(z_0,3R/4)$.  
\qed

\section{ An example}
In this section, we show that if an assumption along the lines of Assumption 2.2(b) does not hold, then  the Harnack inequality can fail.  This example is very similar to the one in \cite{BK1}.  But  since we need some modifications and for the sake of completeness, we give a proof of the following proposition:
\begin{proposition}
There exists a function $n(x,h)$ which satisfies Assumptions 2.2(a) but not (b) and for which the Harnack inequality fails for functions harmonic with respect to the corresponding operator.
\end{proposition}
\proof{
Let $B=B(0,1)$, let $y_0=(1/8,0)$ and for $m\geq4$. let $x_m=(-1/8,2^{-m})$, $z_m=(16, 2^{-m})$, $C_m=B(x_m,2^{-m-4})$, and $E_m=B(z_m,2^{-m-4})$. Define
\[ n(x,h)=\sum_{m=4}^\infty 1_{C_m}(x)1_{E_m}(x+h).\]
Note that $n(x,h)$ satisfies Assumption 2.2(a) and not 2.2(b).  Now we show that $\bP^{y_0}(T_{C_m} < \tau_B)$ is small when $m$ is large.  We see that from Lemma \ref{Lem: exp}, $\bE^{y_0} \tau_B \leq c_1 < \infty$. 
As before we  are going to write $\sL =\sL_c+\sL_d$. Now fix $m$, let $\epsilon = 2^{-m-4}$, let $g(x)=|x-x_m|^{-\beta}$, where $\beta \in (0,1)$ let $\phi$ be a non-negative $C^\infty$ function with support in $B(0,1/2)$ whose integral is 1, let $\phi_\epsilon (x)= \epsilon^{-d}\phi(x/\epsilon)$ and let $f_\epsilon=g * \phi_\epsilon$.  Hence $f_\epsilon \in C^\infty$ and we see that $f_\epsilon \geq c_2 \epsilon^{-\beta}$ on $C_m$.  Since the local part is uniformly elliptic, we have $|\sL_cf_\epsilon(x)|\leq c_3$.  From the definition of $n(x,h)$, we have $|\sL_df_\epsilon(x)|\leq c_4$.  Hence $|\sL f_\epsilon(x)|\leq c_5$.  Since $\bP^{y_0}$ is a solution to the martingale problem for $\sL$, 
\begin{eqnarray*}
\bE^{y_0}f_\epsilon(X_{{T_{C_m}\wedge \tau_B}})-f_\epsilon(y_0)&=&\bE^{y_0} \int_0^{{T_{C_m}\wedge \tau_B}} \sL f_\epsilon(X_s)ds\\
&\leq& c_6\bE^{y_0} \tau_B \leq c_7.
\end{eqnarray*}
Hence
\[c_2\epsilon^{-\beta}\bP^{y_0}(T_{C_m}<\tau_B) \leq \bE^{y_0}f_\epsilon(X_{T_{C_m}\wedge \tau_B})\leq c_7+f_\epsilon(y_0) \leq c_8.\]
Thus $\bP^{y_0}(T_{C_m} < \tau_B)$ will be small if $m$ is large.
Now suppose that the Harnack inequality does hold for non-negative functions that are harmonic in $B$, that is, suppose there exists $c_9$ such that 
\[ u(x)\leq c_9 u(y) \hskip 10mm x,y \in B(0,1/2),\]
for any nonnegative bounded function $u$ which is harmonic in $B$.
Let
\[u_m(x)=\bE^x[1_{E_m}(X_{\tau_B})].\]
Then $u_m(x)$ is bounded. nonnegative, and harmonic in $B$. Note that the only way that $X_{\tau_B}$ can be in $E_m$ is if $X_{\tau_B-}$ is in $C_m$.  We then have, using the assumption that the Harnack inequality holds,
\begin{eqnarray*}
u_m(y_0)&=& \bE^{y_0}[1_{E_m}(X_{\tau_B});T_{C_m}<\tau_B]\\
&=&\bE^{y_0}[\bE^{X_{T_{C_m}}}[1_{E_m}(X_{\tau_B})];T_{C_m}<\tau_B]\\
&=& \bE^{y_0}[u_m(X_{T_{C_m}});T_{C_m}<\tau_B]\\
&\leq& c_9u_m(x_m)\bP^{y_0}(T_{C_m}<\tau_B).
\end{eqnarray*}
Then,
\[\frac{u_m(x_m)}{u_m(y_0)}\geq \frac{1}{c_9\bP^{y_0}(T_{C_m}<\tau_B)}\]
which can be made arbitrary large if we take $m$ large enough.  This is a contradiction and therefore the Harnack inequality cannot hold.
\qed} 

\section*{Acknowledgments}
The author wishes to thank his thesis advisor, Prof. Richard F. Bass, for all his help during the preparation of this work.  The author also thanks an anonymous referee and the associate editor for pointing out several mistakes in earlier drafts of the paper.

\begin{small}
\end{small}


\begin{thebibliography}{99}

\bibitem{Ba2} Bass, R.F (1995).
	\emph{Probabilistic techniques in analysis}. Probability and its Applications. Springer-Verlag, New York.
	%
\bibitem{Ba3} Bass, R.F (1997).
	\emph{Diffusion and Elliptic Operators}. Probability and its Applications. Springer-Verlag, New York.
%
\bibitem{BK1} Bass, R.F and Kassmann, M.(2005). Harnack inequalities for non-local operators of variable order
	\emph{Transactions of the A.M.S}.{\it 357}:837--850
%
\bibitem{BK1} Bass, R.F and Kassmann, M.(2005). Harnack inequalities for non-local operators of variable order
	\emph{Transactions of the A.M.S}.{\it 357}:837--850
%
\bibitem{BK2} Bass, R.F and Kassmann, M.(2005). H\"older continuity of Harmonic functions with respect to operators of variable order.
	\emph{Communications in Part. Diff. Equations}.{\it 7-9}:1249--1259
%
\bibitem{BL1} Bass, R.F and Levin, D.A. (2002). Harnack inequalities for jump processes.
	\emph{Potential Anal.}.{\bf 17}{\it 4}:375--388
%
\bibitem{BL2} Bass, R.F and Levin, D.A. (2002).  Transition probabilities for  symmetric jump processes.
	\emph{Transactions of the A.M.S}.{\bf 357}{\it 7}:2933--2953
%
\bibitem{CK} Chen. Z.Q and Kumagai, T. (2003).  Heat kernel estimates for stable-like processes on d-sets.
	\emph{Stoch. Proc. Their. Applic.,}.{\bf 108}{\it 1}:27--62
%
\bibitem{Fo2}Foondun, M.  Heat kernel estimates and Harnack inequalities for some Dirichlet forms with non-local part.
	\emph{preprint}
%
\bibitem{FS} Fabes. E.B and Stroock, D.W (1986).  A new proof of Moser's Harnack inequality via the old ideas of Nash
	\emph{Arch. Mech. rat. Anal}.{\it 96}:327-338
%	
\bibitem{De} De Giorgi. E (1957).  Sulla differenziabilit{\`a} e l'analiticit{\`a} delle estremali degli integrali multipli regolari
	\emph{Mem. Accad. Sci. Torino. Cl. Sci. Fis. Mat. Nat.}.{\bf 3}{\it 3}:25--43
%
\bibitem{Ja} Jacod. J (1979) .
\emph{Calcul Stochastique et Probl{\`e}mes de Martingales} Lecture notes in Mathematics

%
\bibitem{kas1}Kassmann, M.  The classical Harnack inequality fails for non-local operators.
	\emph{preprint}
%
\bibitem{EL} El Karoui, N and Lepeltier, J.P (1977). Repr\'esentation des processus ponctuels multivari\'es \`a l'aide d'un processus de Poisson
\emph{Z. Wahrscheinlichkeitstheorie und Verw. Gebiete}.{\bf39}{(\it3)}:111--133
%
\bibitem{kry1} Krylov, N.V (1971) An inequality in the theory of stochastic processes.
\emph{Th. Probab. Applic}. {\bf16}: 438--448,
%

\bibitem{KS1}Krylov, N.V  and Safonov, M.V (1979). An estimate for the probability of a diffusion process hitting a set of positive measure
\emph{Dokl. Akad. Nauk SSSR} {\bf 245}{\it 1}:18--20
%
\bibitem{LP}J-P Lepeltier and B. Marchal (1976) Probl{\`e}me des martingales et{ \'e}quations diff{\'e}rentielles stochastiques associ{\'e}es {\`a} un op{\'e}rateur integro-diff{\'e}rentiel
\emph{Ann. Inst H. Poincar{\'e} Sec B (N.S)} {\bf 12}{\it 1}: 43--103
%
\bibitem{Me} Meyer. P.A (1975)Renaissance, recollements, m\'elanges, ralentissement de processus de Markov.
\emph{Ann. Inst. Fourier}{\bf 25 }{\it 3-4}:464--497
%

\bibitem{Mo} Morrey, Ch. B(1938). On the solution of quasi-linear elliptic differential equations
\emph{Trans. Amer. Math. Soc }{\bf }{\it 43 }:126--166
%
\bibitem{M1}Moser. J(1961). On {H}arnack's theorem for elliptic differential equations.
\emph{Comm. Pure Appl. Math. }{\bf }{\it14 }:577--591
%
\bibitem{ M2}Moser. J(1964)A {H}arnack inequality for parabolic differential equations
\emph{ Comm. Pure Appl. Math. }{\bf }{\it 17}:101--134
%
\bibitem{MP}Mikulevicius. R and Pragarauskas. H (1988). On {H}{\"o}lder continuity of solutions of certain integro-differential equations
\emph{Ann. Acad. Scien. Fennicae, Series A. I. Mathematica }{\bf }{\it 13}:  231--238
%
\bibitem{N}Nash, J.(1958)Continuity of solutions of parabolic and elliptic equations.
\emph{Amer. J. Math.  }{\bf 80}{\it }:931--954
%
\bibitem{RRV} Rao. M and Song. R and Vondracek. Z(2006) Green function estimates and {H}arnack inequality for subordinate {B}rownian motions
\emph{Potential Anal. }{\bf 25 }{\it 1}: 1--27
%
\bibitem{St} Stroock. D. W(1975)Diffusion processes associated with {L}\'evy generators.
\emph{Z. Wahrscheinlichkeitstheorie und Verw. Gebiete }{\bf 32}{\it 3}:209--244
%
\bibitem{SV05}R. Song and Z. Vondracek (2005).Harnack inequalities for some discontinuous {M}arkov processes with a diffusion part
\emph{Glas. Math. Ser III. }{\bf 40 }{\it 60}: 177-187.
\end{thebibliography}
\end{document}